\documentclass[journal]{IEEEtran}
\IEEEoverridecommandlockouts
\usepackage{booktabs}
\usepackage{cite}
\usepackage{floatrow}
\usepackage{amsmath,amssymb,amsfonts}
\usepackage{color}
\usepackage{graphicx}
\usepackage{textcomp}
\usepackage{xcolor}
\usepackage{array}

\usepackage[justification=centering]{caption}
\usepackage[ruled,vlined]{algorithm2e}
\usepackage{algpseudocode}
\usepackage{amsmath}
\usepackage{graphics}
\usepackage{epsfig}
\usepackage{caption}
\newsavebox{\ieeealgbox}

\begin{document}

\title{Effective End-to-End Learning Framework for Economic Dispatch}

\author{
Chenbei Lu, Kui Wang \emph{and} Chenye Wu \vspace{-0.0cm}
\thanks{This work has been supported in part by the Zhongguancun Haihua Institute for Frontier Information Technology.}
\thanks{The authors are with the Institute for Interdisciplinary Information Sciences (IIIS), Tsinghua University, Beijing, 100084, P.R. China. C. Wu is the correspondence author. Email: chenyewu@tsinghua.edu.cn}
}

\maketitle

\begin{abstract}
Conventional wisdom to improve the effectiveness of economic dispatch is to design the load forecasting method as accurately as possible. However, this approach can be problematic due to the temporal and spatial correlations between system cost and load prediction errors. This motivates us to adopt the notion of end-to-end machine learning and to propose a task-specific learning criteria to conduct economic dispatch. Specifically, to maximize the data utilization, we design an efficient optimization kernel for the learning process. We provide both theoretical analysis and empirical insights to highlight the effectiveness and efficiency of the proposed learning framework.


\end{abstract}

\label{intro}
\section{Introduction}

The increasing penetration of renewable energies challenges the conventional power system operation paradigm, in particular, the economic dispatch (ED) process. Conventionally, the short term load prediction can be rather accurate, and hence, the dispatch based on those predictions yields near minimal generation cost. However, when renewable generations are considered as negative loads, their stochastic nature makes the net load hard to predict. To improve the effectiveness of ED, various advanced load forecasting methods have been proposed over the past two decades.

Nonetheless, more accurate load forecasting does not necessarily lower the generation cost. This is largely because the common forecasting precision metric is the mean square error (MSE) over the period of interest, and this common metric usually won't coincide with the objective function of ED.

This motivates us to adopt the notion of end-to-end machine learning and propose a task specific criteria to conduct load forecasting. However, to design an effective framework is delicate since end-to-end machine learning usually suffers from low data utilization\cite{b}. To tackle this challenge, we propose an efficient optimization kernel to speed up the training process, which improves the data utilization. The optimization kernel further motivates us to design a robust model-free end-to-end learning framework.

\subsection{Related Works}
We identify two major bodies of closely related research directions in designing the learning framework for ED: one seeks to provide effective load forecasting methods while the other investigates various ways to conduct efficient ED.

Load forecasting is a rather classical technique in power system operation, and has been well investigated (see \cite{r} for a comprehensive survey). While most classical methods for load forecasting utilize the statistical analytics (e.g., the adaptive autoregressive moving-average (ARMA) model in \cite{l1}) or time series analysis (e.g., stochastic time series analysis in \cite{l3}), machine learning algorithms have also been applied to load forecasting since the mid-1980s \cite{lf1}. With the advance in machine learning over the last decade, this line of research attracts significant attention.
Just to name a few, Lee \emph{et al}. design a neural network based model to capture weekend-day energy consumption patterns in \cite{l4}. Bashir \emph{et al}. seek to improve the effectiveness of neural network by Particle Swarm Optimization (PSO) in \cite{l5}.
Javed \emph{et al.} combine the Artificial Neural Network (ANN) and Support Vector Machine (SVM) for load forecasting in \cite{l7}.
We want to emphasize that the conventional forecasting precision metric in the literature is MSE. In our work, we aim to highlight the mismatch between MSE and the desirable task specific criteria.


Economic dispatch is one of the most important processes in power system operation. There is a huge body of related literature to improve the effectiveness of ED. The major difficulties come from the temporal coupling and the dynamic implementation.
To tackle temporal challenges, the solution concepts range from classical linear and quadratic programming \cite{b3} to genetic algorithm for value point loading \cite{b4}. To overcome the difficulties of dynamic dispatch, various mathematical programming approaches have been proposed, including Lambda iterative method \cite{b7}, interior point method \cite{b8,b9}, and dynamic programming \cite{b10}. However, mathematical programming approaches are usually time consuming and not ideal for large scale systems. Hence, for large scale non-convex ED, heuristic and hybrid methods are often preferable, e.g., hybrid evolution programming in \cite{b11}.

However, the notion of using learning framework to conduct ED only appears recently.
Donti \emph{et al}. propose a generic end-to-end machine learning framework for stochastic programming, with an application to the single generator ED problem in \cite{z}. Different from this work, we target to design an end-to-end machine learning framework for general ED.

\subsection{Our Contributions}
In seek of designing the effective end-to-end machine learning framework, our principal contributions can be summarized as follows.


\begin{itemize}
    \item \emph{Task Specific Criteria}:  We identify that MSE is not ideal to evaluate the performance of load forecasting for conducting ED. We adopt the notion of end-to-end machine learning to derive the task specific criteria.
    \item \emph{Optimization Kernel Construction}: To avoid solving the multivariate constrained stochastic optimization problem during each learning iteration, we exploit the problem structure and propose an efficient way to construct the optimization kernel for effective learning.
    {
    \item \emph{Model-Free End-to-End Learning}: Motivated by the optimization kernel, we propose a model-free approach to further improve the efficiency and effectiveness of our end-to-end learning framework.}
\end{itemize}

The rest of the paper is organized as follows. 
Section \ref{formulation} presents the system model. We revisit the conventional learning framework for ED in Section \ref{conventional_frame}. Then, in Section \ref{e2elearning_basis}, we lay out the theoretical foundation for the end-to-end learning framework for ED, which can be time consuming. To tackle the challenge, we design the optimization kernel in Section \ref{implementation}, and then propose the model-free framework in Section \ref{modelfree}.
Numerical studies verify the efficiency and effectiveness of our proposed frameworks in Section \ref{experiment}. Finally, concluding remarks and future directions are given in Section \ref{conclude}. We defer all the necessary proofs in the Appendix.

\section{Problem Formulation}
\label{formulation}
We consider a standard ED process over a period of $T$ time slots with geographically distributed $n$ generators and $m$ demands. Besides the generation cost for each generator, we assume there are certain risk costs associated with the supply demand mismatch in real time. To capture the stochastic nature of net demands, they are modeled as random variables, following possibly distinct distribution. The system operator seeks to minimize the total cost of the system by solving the following optimization problem :
\begin{align}
    \min \quad &\sum_{t=1}^{T}{\left(\sum_{i=1}^{n}c_i({g}_{it})+\gamma_1{\mathbb{E}(S_t)}^++\gamma_2{\mathbb{E}(-S_t)}^+\right)}
    \label{eq1}\\
    s.t. \quad 
    &0 \leq g_{it} \leq B_i, \quad \forall{i,t}
    \label{eq2}\\
    &-\boldsymbol{b} \leq \mathbb{E}(H_g\boldsymbol{g}_{t}-H_d\boldsymbol{\hat{\boldsymbol{d}}_t}) \leq \boldsymbol{b}, \quad
    \forall{t}
    \label{eq3}\\
    &S_{t} = \sum \nolimits_{j=1}^m{\hat{d}_{jt}}-\sum\nolimits_{i=1}^n{g_{it}},
     \quad \forall{t}.
     \label{eq4}
\end{align}

The decision variables are $g_{it}$'s, and each $g_{it}$ denotes the dispatched generation of generator $i$ at time $t$.  We denote $\boldsymbol{g}_t$ the vector of $(g_{1t},...,g_{nt})$. On the other hand, the parameters in the optimization problem are defined as follows.
\begin{itemize}
    \item $\hat{d}_{jt}$ : prediction of demand $j$ at time $t$, $\hat{\boldsymbol{d}_t} = (\hat{d}_{1t},...,\hat{d}_{mt})$.
    \item $\gamma_1$, $\gamma_2$ : unit generation shortage and excess penalties
    \item $S_t$ : total shortage based on predicted demands at time $t$
    \item $B_i$ : generation capacity of generator $i$
    \item $\boldsymbol{b}$ : transmission line capacity vector
    \item $H_g$, $H_d$: generation and load shift factor matrices
    \item $c_i(\cdot)$ : generation cost function of generator $i$
\end{itemize}

Constraint (\ref{eq2}) captures the limited capacity of each generator. Constraint (\ref{eq3}) uses the shift factor matrices \cite{mp} to represent the DC approximation of line capacity constraints.
We follow the literature and assume the generation cost functions (i.e., $c_i(\cdot)$'s) are linear for each individual generator. Hence, the total individual cost function is piecewise linear.
Note that, the randomness in the predicted demands ($\hat{d}_{jt}$'s) makes the optimization problem (\ref{eq1})-(\ref{eq4}) a multivariate stochastic optimization with linear constraints.

\section{Conventional Learning Framework for ED}
\label{conventional_frame}
In this section, we first introduce the conventional learning framework for ED and then use the electricity pool model to highlight its drawback: the mismatch between MSE learning criteria and the ultimate goal of minimizing system cost.
\subsection{Conventional Framework}
\begin{figure}[t]
    \centerline{\includegraphics[width=8cm]{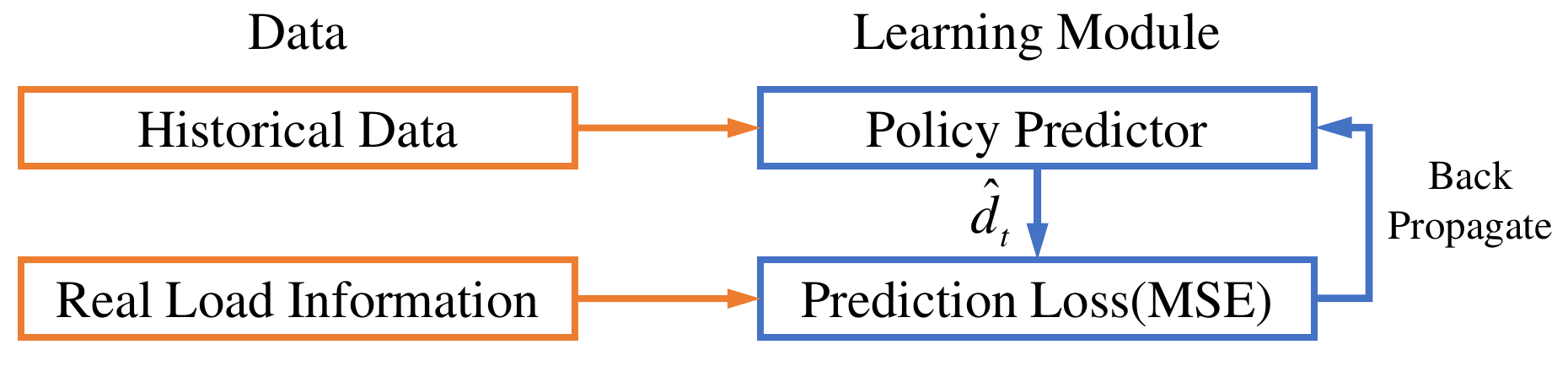}}
    \centering
\caption{Conventional Learning Process.}
\label{fig1}
\end{figure}

Conventional wisdom isolates the load prediction from the whole ED process and aims to obtain a perfect load predictor across time. Hence, the major task of conventional learning framework for ED is to train an accurate predictor. And the ED process will directly take the predicted loads as single value inputs and conduct the dispatch. The learning process is visualized in Fig.\ref{fig1}.
The whole process is intuitive and easy to understand. 
However, the training criteria (MSE) is a generic selection, which is not customized for ED. Hence, the load predictor trained through this process does not necessarily lower the generation cost in practice.



\subsection{Inefficiency of MSE based Predictor}

To highlight the fact that MSE can be inefficient, we consider the ED in the electricity pool model, where all the network constraints are ignored and the system cost function degenerates to a single piecewise linear function.

Mathematically, we simplify the optimization problem (\ref{eq1})-(\ref{eq4}) as follows:
\begin{equation}
    \begin{split}
         \min \quad &\sum\nolimits_{t=1}^{T}\left[\hat{C}_t({g}_{t})+\gamma_1{\mathbb{E}(S_t)}^++\gamma_2{\mathbb{E}(-S_t)}^+\right]\\
    s.t. \quad&S_{t} = \sum \nolimits_{j=1}^m{\hat{d}_{jt}}-\sum\nolimits_{i=1}^n{g_{it}},
     \quad \forall{t},
    \end{split}
    \label{elecpool}
\end{equation}
where for each $t$,
\begin{equation}
    \begin{split}
        \hat{C}_t(g_t) = \min\quad&\sum\nolimits_{i=1}^n{c_i(g_{it})}\\
    s.t.\quad&0 \leq g_{t} := \sum\nolimits_{i=1}^n{g_{it}} \leq \sum\nolimits_{i=1}^n B_i.
    \end{split}
\end{equation}

Note that, in optimization problem (\ref{elecpool}), we only need to consider the total dispatched generation $g_t$'s, instead of $g_{it}$'s, as it is easy to recover the vector $(g_{it}, \forall{i})$ from $g_t$ based on the merit order of the generator's marginal cost. It is also clear that in this simplified model, the dispatch decisions $g_t$'s are decoupled across time. To capture the role of load predictor in the decision making, we define
\begin{align}
    \hat{d}_t = \sum\nolimits_{j=1}^m{\hat{d}_{jt}} . 
\end{align}

We further denote $f_t(x)$ and $F_t(x)$ the probability density function (\emph{pdf}) and the cumulative density function (\emph{cdf}) of $\hat{d}_t$. They allow us to express the risk cost $\hat{R}(g_t,f_t(x))$ as follows.
\begin{equation}
\begin{aligned}
    \hat{R}(g_t,f_t(x))=& \gamma_1\int_{g_t}^{\infty}{(x-g_t)f_t(x)dx} \\
    &+ \gamma_2\int_{-\infty}^{g_t}{(g_t-x)f_t(x)dx}.
    \end{aligned}
\end{equation}
Based on $\hat{R}(g_t,f_t(x))$, we can rewrite the decision making problem for each time $t$:
\begin{equation}
    \begin{split}
        \min \quad & \hat{C}_t({g}_{t})+ \hat{R}(g_t,f_t(x)) \\
    s.t. \quad &0 \leq g_{t} \leq \sum\nolimits_{i=1}^n B_i.
    \end{split}
\end{equation}

We assume the generation shortage penalty is larger than the marginal generation cost, and far overweigh the excessive penalty, i.e., $\gamma_1 \!>\! C'(g_t)$, $\gamma_1 \gg \gamma_2$.
This is reasonable as $\gamma_1$ represents the marginal generation cost for the next more expensive generator while $\gamma_2$ can be understood as the opportunity cost for the undispatched generators.
The first order optimality condition indicates the optimal dispatch should satisfy the following condition:
\begin{align}
     g_t = \min \left\{F_t^{-1}\left(\frac{\gamma_1-\hat{C}_t'(g_t)}{\gamma_1+\gamma_2}\right), \sum\nolimits_{i=1}^n{B_i}\right\}.
\end{align}

Under normal operating conditions, the available generation ($\sum\nolimits_{i=1}^n{B_i}$) is always greater than the peak demand. Hence, we can safely drop the min operator in practice:
\begin{align}
    g_t = F_t^{-1}\left(\frac{\gamma_1-\hat{C}_t'(g_t)}{\gamma_1+\gamma_2}\right).
\end{align}

Following the notion that generation and the predicted demand should be aligned, this notion indicates that if we are only allowed to obtain a single valued load prediction, we should deliver $d^*_t$ to the system operator, where
\begin{align}\label{desired_criteria}
    d_t^* = F_t^{-1}\left(\frac{\gamma_1-\hat{C}_t'(d_t^*)}{\gamma_1+\gamma_2}\right).
\end{align}

On the other hand, the learning framework based on MSE criteria will deliver $d_t^* = \mathbb{E}[\hat{d}_t]$ to the system operator, which is often different from the percentile\footnote{Although some precision metrics like mean absolute error (MAE) will lead to a percentile form for $d^*_t$, they are often independent of $\gamma_1$, $\gamma_2$, and the cost functions. Hence, they are also fundamentally different from task specific criteria.} represented in Eq. (\ref{desired_criteria}). Hence, the conventional framework is not efficient in minimizing the system cost even in the simplified electricity pool model.

\section{End-to-End Learning : the Basis}
\label{e2elearning_basis}



End-to-end learning is a powerful tool in the machine learning community. In this section, we revisit the basic concepts in the generic end-to-end learning framework, and then introduce the notion of task specific criteria for ED.

Specifically, as for learning the ED policies, end-to-end learning would directly learn the final dispatch policies given the input data and the loss function for training, which often measures the difference between predicted policy $\hat{\boldsymbol{g}}_{t}$ and the optimal policy ${\boldsymbol{g}}^*_{t}$, i.e, $\sum\nolimits_{t = 1}^{T}\sum\nolimits_{i = 1}^n|\hat{{g}}_{it} - g^*_{it}|$.

However, ignoring all the intermediate stages makes pure end-to-end learning suffer from a number of disadvantages.
First, the predicted dispatch policy may not satisfy all the constraints in ED. Also, this framework can be rather data inefficient, especially for large dynamic system control, such as ED in power system.

\subsection{Task Specific Criteria for ED}
In order to avoid the shortcomings of pure end-to-end learning, 
a task specific criteria can be helpful.
Specifically,
we can carefully examine the structure of problem (\ref{eq1})-(\ref{eq4}).
Given the probability distribution of all the random variables, we can convert this stochastic optimization into a deterministic optimization.

Note that, since our formulation ignores the ramping constraints, the decision makings across time are naturally decoupled. Such simplifications allow us to highlight the effectiveness of task specific end-to-end machine learning. 
Hence, in the subsequent analysis,
we will focus on single shot ED and provide our insights on how to generalize our approach to ED with ramping constraints by the end of this section. Suppose we are given the distribution of $\hat{d}_t$, then the single time shot ED degenerates to the following problem:
\begin{equation}
\begin{split}
    \min \quad &\sum_{i=1}^{n}c_i({g}_{it})+\gamma_1{\mathbb{E}(\hat{d}_t-g_t)}^++\gamma_2{\mathbb{E}(g_t-\hat{d}_t)}^+\\
    s.t. \quad 
    &0 \leq g_{it} \leq B_i, \quad \forall{i}\\
    &-\boldsymbol{b} \leq \mathbb{E}(H_g\boldsymbol{g}_{t}-H_d\boldsymbol{\hat{\boldsymbol{d}}_t}) \leq \boldsymbol{b}. 
\end{split}
\label{ED_tsk}
\end{equation}

\vspace{0.1cm}
\noindent \textbf{Lemma 1}:
The single shot ED problem (\ref{ED_tsk}) is convex.
\vspace{0.1cm}

The proof can be immediately obtained by examining the Hessian matrix of the objective function with respect to ${g_{it}}$'s and showing it is positive semi-definite.
The objective function in (\ref{ED_tsk}) can be used to derive the chain rule for back propagation in the task specific end-to-end learning.

Suppose $\hat{d}_t$ follows some distribution $f_t(x;\boldsymbol{\hat{\theta}}_t)$, where $\boldsymbol{\hat{\theta}}_t$ represents the parameters of the distribution (e.g., parameter $\lambda$ for Exponential distribution). Then, the goal of the task specific predictor is to estimate the distribution $f_t$ with the appropriate parameters $\boldsymbol{\hat{\theta}}_t$. Given the predicted distribution, we define the task loss function by $L_t(\boldsymbol{\hat{g}}_t,d_t)$:
\begin{equation}
\!L_t(\boldsymbol{\hat{g}}_t,d_t) \!=\! C(\boldsymbol{\hat{g}}_t)\!+\gamma_1{({d}_t-\hat{g}_t)^+}\!+\gamma_2(\hat{g}_t-{d}_t)^+\!-C_t(d_t),
\end{equation}
where $d_t$ denotes the true load.
The other notations $C(\boldsymbol{g}_t)$, $C_t(d_t)$, and $\hat{\boldsymbol{g}}_t$ represent the total generation cost function given dispatch profile $\boldsymbol{g}_t$, the minimal generation cost to meet the true demand $d_t$, and the optimal generation dispatch profile given estimation $f_t(x;\hat{\boldsymbol{\theta}}_t)$. Formally,
\begin{align}
    C({\boldsymbol{g}}_t) = \sum\nolimits_{i = 1}^n c_i({g}_{it}),
\end{align}
\begin{equation}
    \begin{split}
        C_t(d_t) := \min \quad &\sum\nolimits_{i=1}^{n}c_i({g}_{it})\\
    s.t. \quad&\text{Constraints }(\ref{eq2}),(\ref{eq3})\\
     &\sum\nolimits_{i=1}^n g_{it} = d_t,
    \end{split}
    \label{optkernel}
\end{equation}
\begin{equation}
    \begin{split}
    \hat{\boldsymbol{g}}_t = \arg\min \quad &C(\boldsymbol{g}_t) +\hat{R}(g_t,f_t(x;\hat{\boldsymbol{\theta}}_t))\\
    s.t. \quad &\text{Constraints }(\ref{eq2}),(\ref{eq3}).
    \end{split}
\end{equation}

Note that, $\hat{\boldsymbol{g}}_t$ can be calculated by sequential quadratic programming (SQP). Hence, the chain rule of back propagation aligns with the chain rule for partial derivatives:
\begin{align}
    \frac{\partial{L}}{\partial{\hat{\boldsymbol{\theta}}_t}} = \frac{\partial{L_t}}{\partial{\hat{\boldsymbol{g}}_t}}
    \cdot\frac{\partial{\hat{\boldsymbol{g}}_t}}{\partial\hat{\boldsymbol{\theta}}_t} , \forall{t}.
    \label{chain}
\end{align}

\begin{figure}[t]
\centerline{\includegraphics[width=9cm]{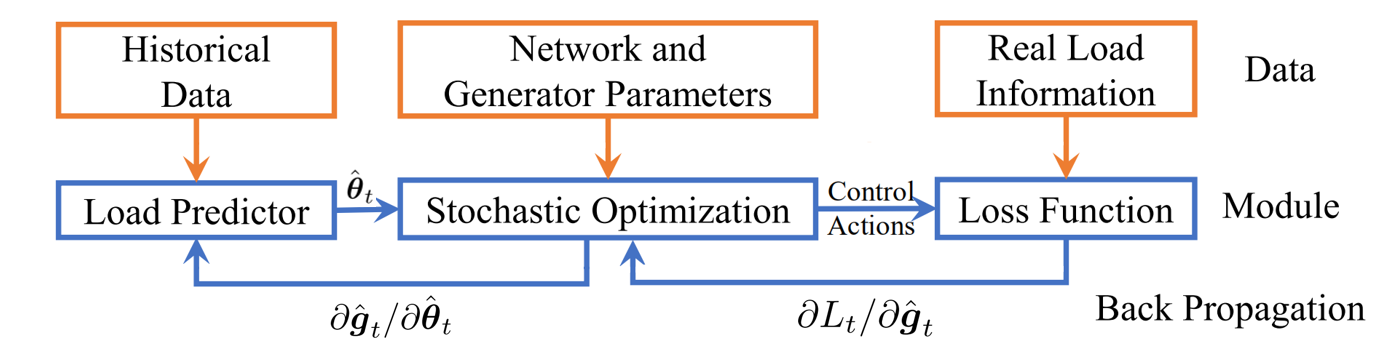}}
\caption{Task-Specific Optimization Based Learning Process.}
\label{e2e4p}
\end{figure}

The two partial derivatives on the right hand side in Eq.(\ref{chain}) can be calculated automatically in the SQP (see \cite{opt} for more details) without explicit formula.
We visualize the whole process in Fig. \ref{e2e4p}.

\section{End-to-End Learning : Implementation}
\label{implementation}
In this section, we discuss in detail how to implement the task specific criteria to achieve the best efficiency and effectiveness. Specifically, by observing the structure of the single shot ED problem, we design an optimization kernel to enable efficient learning.

\subsection{Optimization Kernel for Efficient Learning}

We first use the notion of parametric analysis to examine problem (\ref{ED_tsk}).
The key parameter in determining the value of the objective function is the same as that for the electricity pool model, the total generation $g_t$. We define $\tilde{L}_t(g_t)$ to highlight this observation:
{
\begin{equation}
\begin{split}
        \tilde{L}_t(g_t)\! :=\! \min\quad&C(\boldsymbol{g}_t)\!+\!\gamma_1{\mathbb{E}(\hat{d}_t\!-\!g_t)}^+\!+\!\gamma_2{\mathbb{E}(g_t\!-\!\hat{d}_t)}^+\\
        s.t.\quad&\text{Constraints }(\ref{eq2}),(\ref{eq3}),\\
    &g_t = \sum\nolimits_{i = 1}^{n}g_{it}.
\end{split}
\label{Lgt}
\end{equation}
}

$\tilde{L}_t(g_t)$ can be decomposed into two components. One represents the cost of a fixed structure (the generation cost, $C_t(g_t)$ as defined in Eq.(\ref{optkernel})) and the other deals with the predicted distribution of $\hat{d}_t$ (the risk cost, denoted by $\tilde{R}_t(g_t)$). 
Mathematically, 
\begin{align}
    \tilde{L}_t(g_t) = C_t(g_t) + \tilde{R}_t(g_t),
\end{align}
where
\begin{align}
\tilde{R}_t(g_t) = &\gamma_1\mathbb{E}{(\hat{d_t}-g_t)^+} + \gamma_2\mathbb{E}{(g_t-\hat{d_t})^+}.   
\end{align}

Denote the feasible region constructed by constraints of (\ref{Lgt}) by $\mathcal{A}(g_t)$.
Define
\begin{align}
&g_t^{\min} = \inf \mathcal{A}(g_t),\\
&g_t^{\max} = \sup \mathcal{A}(g_t).
\end{align}

The following lemma indicates the structure of $C_t(g_t)$.

\vspace{0.1cm}
\noindent \textbf{Lemma 2}: The generation function $C_t(g_t)$ is continuous, piecewise linear, and convex in $g_t$ over $[g_t^{\min}, g_t^{\max}]$.
\vspace{0.1cm}

Inspired by the idea proposed in \cite{w}, we can efficiently construct $C_t(g_t)$ as described in Algorithm \ref{CGAs} (Curve Generative Algorithm or CGA in short). 

\begin{algorithm}
\SetAlgoLined
\SetKwInOut{Input}{Input}
\SetKwInOut{Output}{Output}

\Input{generation range: $g_t^x$ and $g_t^y$\;}
\Output{function curve of $C_t(g_t)$ over $[g_t^x, g_t^y]$\;}
Solve $C_t(g_t^x)$, $C_t(g_t^y)$ and get their Lagrangian multipliers $\lambda_t^x$, $\lambda_t^y$, respectively;

\eIf {$\lambda_t^x$ == $\lambda_t^y$}
{The slope of $C_t(g_t)$ over $[g_t^x, g_t^y]$ is $\lambda_t^x$;}
{Solve the equations below to obtain $g_t^z$, $C_z$;
\begin{equation}
\left\{
\begin{aligned}
&C_z - C_t(g_t^x) = \lambda_t^x(g_t^z-g_t^x)\\
&C_z - C_t(g_t^y) = \lambda_t^y(g_t^z-g_t^y)
\end{aligned}
\right.
\end{equation}
}
\eIf {$C_z$==$C_t(g_t^z)$}
{$g_t^z$ is the unique breaking point in $[g_t^x, g_t^y]$, the
slope over $[g_t^x, g_t^z]$ is $\lambda_t^x$ and the slope of $[g_t^z, g_t^y]$ is $\lambda_t^y$;}
{Recursively generate the curve by \text{CGA}$(g_t^x,g_t^z)$ and \text{CGA}$(g_t^z,g_t^y)$;}
Return;
\caption{CGA($g_t^x$, $g_t^y$)}
\label{CGAs}
\end{algorithm}


With the knowledge of $C_t(g_t)$, the optimal total dispatch $g_t^*$ should satisfy:
\begin{equation}
\begin{split}
    g^*_t = \arg\min&\quad \tilde{L}_t(g_t)\\
    s.t.&\quad g_t^{\min} \leq g_t \leq g_t^{\max}.
\end{split}
\label{g^*_t}
\end{equation}
It is straightforward to verify that
\begin{align}
    C_t''(g_t) + \tilde{R}_t''(g_t) > 0, \quad g_t \geq 0.
\end{align}

{Hence, there must exist a unique $g^*_t$ to problem (\ref{g^*_t}).}
While $C_t(g_t)$ is not differentiable everywhere over the domain, we can conduct a binary search to efficiently obtain $g^*_t$ for given prediction distributions.

The only remaining handle is to construct the dispatch policy $\boldsymbol{g}_t$ from the total dispatch ${g}_t$. Due to the continuity of the solution space of $C_t(g_t)$, we can efficiently construct $\boldsymbol{g}_t$ with the help of the information embedded in the breaking points obtained in Algorithm \ref{CGAs}. More precisely, without loss of generality, suppose $g^*_t \in [g_t^k,g_t^{k+1}]$, where $g_t^k$ and $g_t^{k+1}$ are two adjacent breaking points of $C_t(g_t)$. Denote $\boldsymbol{g}_t^k$ and $\boldsymbol{g}_t^{k+1}$ the corresponding dispatch profiles for ${g}_t^k$ and ${g}_t^{k+1}$. Then, continuity of solution space leads to the following lemma:

\vspace{0.1cm}
\noindent \textbf{Lemma 3}:
If $g^*_t = (1-\gamma)g^k_t + \gamma{g^{k+1}_t}, 0 \leq \gamma \leq 1$, then the optimal dispatch profile $\boldsymbol{g}_t^*$ corresponding to the total dispatch $g^*_t$ can be constructed as follows:
\begin{align}
    \boldsymbol{g}_t^* = (1-\gamma)\boldsymbol{g}^k_t + \gamma{\boldsymbol{g}^{k+1}_t}.
\end{align}
\vspace{0.1cm}

Thus, we complete the construction of the optimization kernel for the end-to-end framework. We illustrate this process in detail in Algorithm \ref{kernel}, where we utilize the information of $C'_t(g_t)$. This information can be obtained during the construction process of $C_t(g_t)$. Since $C'(g_t)$ is only well defined over $[g_t^{\min}, g_t^{\max}]$, we generalize this derivative beyond $[g_t^{\min}, g_t^{\max}]$ as follows:
\begin{align}
  C_t'(g_t)=\left\{
\begin{aligned}
&C_t'(g_t^{\min}),   &g_t \leq g_t^{\min} \\
&C_t'(g_t),   &g_t^{\min} < g_t  < g_t^{\max} \\
&C_t'(g_t^{\max}),   &g_t \geq g_t^{\max}
\end{aligned}
\right.  
\end{align}
\footnotetext[2]{Each piece of data in the two sets can be represented by a pair $(x,d)$, where $x$ denotes the input data, like the previous day's load, week-weekend types and temperatures, while $d$ denote the output data, like the true load of the next day.}
\begin{algorithm}
\SetAlgoLined
\SetKwInOut{Input}{Input}
\SetKwInOut{Output}{Output}

\Input{training set and validation set $(x,d)$\footnote[2];}
\Output{predictor $P$;}
Compute $C_t'(g_t)$'s and $\tilde{R}_t'(g_t)$'s function curve;

Compute $g^{\min} $ and $g^{\max} $;

Compute ${G}(g_t)$'s function curve, where $\boldsymbol{g}_t = G(g_t)$;

\While{Loss of validation set $L_v$ doesn't increase}
{
Sample a batch of $(x,d)$ from training set;\\
Estimate $\hat{\boldsymbol{\theta}}_t$ from $P$ and input data $x$:
$$\hat{\boldsymbol{\theta}_t} = P(x)$$\\
Generate distribution $f_t(\hat{\boldsymbol{\theta}}_t)$;\\
$\hat{d}_t = g^{\min} $;\\
$gap = g^{\min} - g^{\max}  $;\\ 
\While {$gap > (g^{\min}-g^{\max})\cdot10^{-6}$ }
{
$\hat{d}_t = \hat{d}_t - sign(C'(\hat{d}_t)+\tilde{R}'(\hat{d}_t))\cdot{gap}$;\\
$gap = gap/2$;
}
$\hat{g}_t = \text{med}\left\{g^{\min}, \hat{d}_t, g^{\max}\right\}$;\\
Compute generation policy $\hat{\boldsymbol{g}}_t$ from ${G}(\hat{g}_t)$;\\
Compute task loss $L_t(\hat{\boldsymbol{g}}_t,d)$\\
Compute derivatives $\frac{\partial{L_t(\hat{\boldsymbol{g}}_t,d)}}{\partial\hat{{\boldsymbol{\theta}}}_t}$;\\
Back propagate and update the predictor $P$;\\
Compute the loss in validation set $L_v$ by $P$;
}
\caption{Optimization Kernel Based  Learning}
\label{kernel} 
\end{algorithm}

We want to close this section by making remarks on generalizing our framework to the ED problem with ramping constraints. The only difference is that the dispatch policies cannot be decoupled across time. Hence, the risk cost has to rely on $T$ variables, i.e., $g_1, ..., g_T$, instead of a single variable. In this case, we need to construct a high dimensional optimization kernel to enable the efficient learning.

\section{End-to-End Learning: Model-Free}
\label{modelfree}
While the optimization kernel can effectively speed up the learning process, in this section, we highlight that the knowledge of distribution is great, but not essential.

Denote $f(d_t;\hat{\boldsymbol{\theta}}_t)$ and $F(d_t;\hat{\boldsymbol{\theta}}_t)$ the \emph{pdf} and \emph{cdf} of hypothetical distribution of $\hat{d}_t$, and denote $H(d_t)$ its true \emph{cdf}.

Thus, the derivative of task loss ${L}_t(g_t)$ with respect to $g_t$ can be obtained as follows:
\begin{align}
    \frac{\partial{{L}_t(g_t)}}{\partial{g_t}} = (\gamma_1 + \gamma_2)F(g_t;\hat{\boldsymbol{\theta}}_t) - \gamma_1 + C'(g_t).
\end{align}
Combining the first order optimality condition with the generation capacity constraint $g_t \in [g_t^{\min}, g_t^{\max}]$, the optimal policy $g^0_t$ given the distribution $F(x;\hat{\boldsymbol{\theta}}_t)$ can be obtained as follows.
\begin{align}
    g_t^0 &= \text{med}\left\{g_t^{\min},F^{-1}\left(\frac{\gamma_1-C'(g_t^0)}{\gamma_1+\gamma_2};\hat{\boldsymbol{\theta}}_t\right), g_t^{\max}\right\}.
    \label{g0}
\end{align}
Note that $\text{med}$ is median operator. We denote $g_t^0 = g_t(\hat{\boldsymbol{\theta}}_t)$ to highlight that it could be a function of the distribution parameters $\hat{\boldsymbol{\theta}}_t$. The same analysis applies to decide the true optimal control policy $g_t^*$ (given $H(d_t)$):
\begin{align}
    g_t^* &= \text{med}\left\{g_t^{\min},H^{-1}\left(\frac{\gamma_1-C'(g_t^*)}{\gamma_1+\gamma_2}\right), g_t^{\max}\right\}.
    \label{g*}
\end{align}
Combining Eqs. (\ref{g0}) and (\ref{g*}), we can examine the performance loss $\Delta_L(g_t^*, g_t(\hat{\boldsymbol{\theta}}_t))$ induced by the inaccurate distribution estimation:
\begin{equation}
    \begin{split}
        \Delta_L &= |{L}_t(\hat{\boldsymbol{\theta}}_t)-{L}_t(g_t^*)|= \left|\int_{g_t^*}^{g_t(\hat{\boldsymbol{\theta}}_t)}{L_t'(x)dx}\right|\\
    &= \left|\int_{g_t^*}^{g_t(\hat{\boldsymbol{\theta}}_t)}{(\gamma_1+\gamma_2)H(x)-\gamma_1+C_t'(x)dx}\right|.
    \end{split}
\end{equation}
Define 
\begin{align}
   K(x) = (\gamma_1 + \gamma_2)H(x) - \gamma_1 + C_t'(x),
\end{align}
it is straightforward to verify the following lemma:

\vspace{0.1cm}
\noindent \textbf{Lemma 4}:
The function $K(x)$ is either always positive over $[g_t^* , g_t(\hat{\boldsymbol{\theta}}_t)]$ or always negative over $[g_t(\hat{\boldsymbol{\theta}}_t), g_t^*]$.
\vspace{0.1cm}

With the lemma, we can conclude that the ultimate goal of the learning process should be to estimate $g^*_t$ as accurately as possible, which does not require the explicit knowledge of distribution. This is also desirable according to Occam's Razor \cite{Occam}: estimating the distribution often leads to more parameter estimations.
Specifically, when the data are not sufficiently large, the trained model may suffer from overfitting issues or poor generalization ability in practice.

To achieve an accurate prediction for $g^*_t$, we design the following loss function for the learning process:
\begin{align}
    Q_t(\hat{g}_t, d_t) \!=\! C_t(\hat{g}_t)\!-\! C_t(d_t)\!+\!\gamma_1(d_t\!-\! \hat{g}_t)^+ \!+\! \gamma_2(\hat{g}_t- d_t)^+.
    \label{qt}
\end{align}
We want to emphasize that in Eq.(\ref{qt}), $d_t$ is the actual demand during learning.
Also, the geographical distributed demand information has been encoded into the function $C_t(\cdot)$ via $\mathbb{E}[d_{jt}]$, $\forall{j}$, together with all the transmission line capacity constraints. Based on this loss function, we design the model-free end-to-end learning framework as described in Algorithm \ref{mdfree}.
\begin{algorithm}
\SetAlgoLined
\SetKwInOut{Input}{Input}
\SetKwInOut{Output}{Output}

\Input{training set and validation set $(x,d)$;}
\Output{predictor $P$;}
Compute $C_t(g_t)$'s function curve;

Compute $g^{\min} $ and $g^{\max} $;

\While{Loss in validation set $L_v$ doesn't increase}
{
Sample a batch of $(x,d)$ from training set;\\
Compute the total generation $\hat{g_t}$ from predictor $P$ and input data $x$:
$$\hat{g_t} = P(x)$$\\
$\hat{g_t} = \text{med}\left\{g^{\min}, \hat{g_t}, g^{\max}\right\}$;\\
Compute task loss $Q_t(\hat{g_t},d)$;\\
Compute derivatives $\frac{\partial{Q_t(\hat{{g}}_t,d)}}{\partial{\hat{{g}}_t}}$;\\
Back propagate and update the predictor $P$;\\
Compute the loss in validation set $L_v$ by $P$\;
}
\caption{Optimization-Free Learning}
\label{mdfree} 
\end{algorithm}

\section{Numerical Studies}
\label{experiment}
We compare the performance of different approaches on a four-bus system in terms of effectiveness, robustness and efficiency. We further illustrate the effectiveness of the model-free framework on the IEEE 39-bus system.
\begin{figure}[t]
\centerline{\includegraphics[width=4cm]{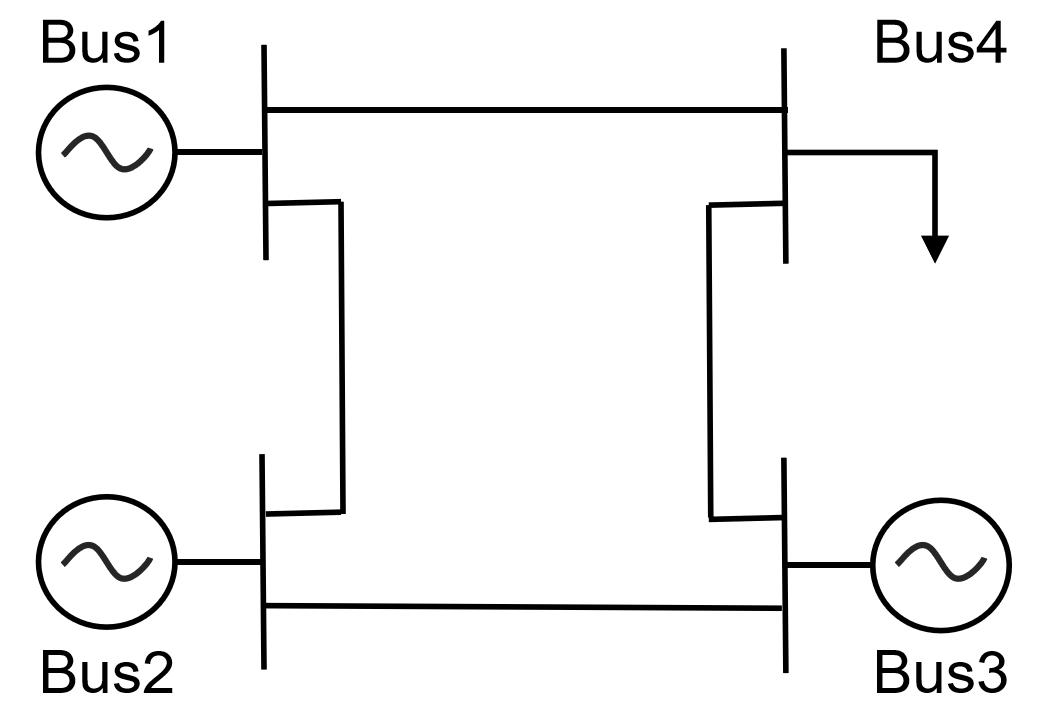}}
\caption{4-Bus System.}
\label{4bus}
\end{figure}
\begin{table}[t]
	\centering
	\caption{4-Bus System Parameters}
	\begin{tabular}{ccccc}
		\toprule  
		Bus-Bus  & 1 & 2  &  3 &  4 \\ 
		\midrule  
		Impedance (p.u.) & $-j$ & $-j$ & $-j$ & $-j$ \\
		\midrule  
		Capacity (MW) & 1.5 & 1.5 & 1.5 & 1.5\\
		\bottomrule  
	\end{tabular}
	\label{4bus_grid}
\end{table}
\subsection{Setup for 4-Bus System}
As shown in Fig. \ref{4bus}, the system contains 3 generators at Bus 1 to Bus 3, respectively, and a single load at Bus 4. To capture the stochastic nature of the net demand, we use the 5-year PJM load data from 2012 to 2016 \cite{PJM}. We assume the marginal costs of the three generators are $\$$40, $\$$50, $\$$60/MWh, respectively.
We further set the unit shortage penalty $\gamma_1$ to be $\$$100/MWh, and the unit excess penalty $\gamma_2$ to be $\$$10/MWh. The impedances and the line capacities of the network can be found in Table \ref{4bus_grid}.

We employ a 2-layer neural network with 128 neurons in each layer to implement our proposed frameworks. During the learning process, we divide the 5-year data set into 3 sets: training set, validation set, and test set. The training set contains the load data in the first $1,200$ days. The data in the following 200 days construct the validation set, and the data in the last 400 days form the test set.

All models are trained to the convergence of accuracy decided by the validation set.
For all the models, they share the same inputs: historical hourly load data in the former 24 hours and weekday-weekend type. 



\subsection{Performance Evaluation: Effectiveness}
We compare the performance of four frameworks: conventional approach with MSE criteria, end-to-end learning with task specific criteria, end-to-end learning with optimization kernel, and model-free end-to-end framework. The evaluation metrics are load prediction error (MSE) and the loss in dispatch cost.

\begin{figure}[t]
\centerline{\includegraphics[width=8cm]{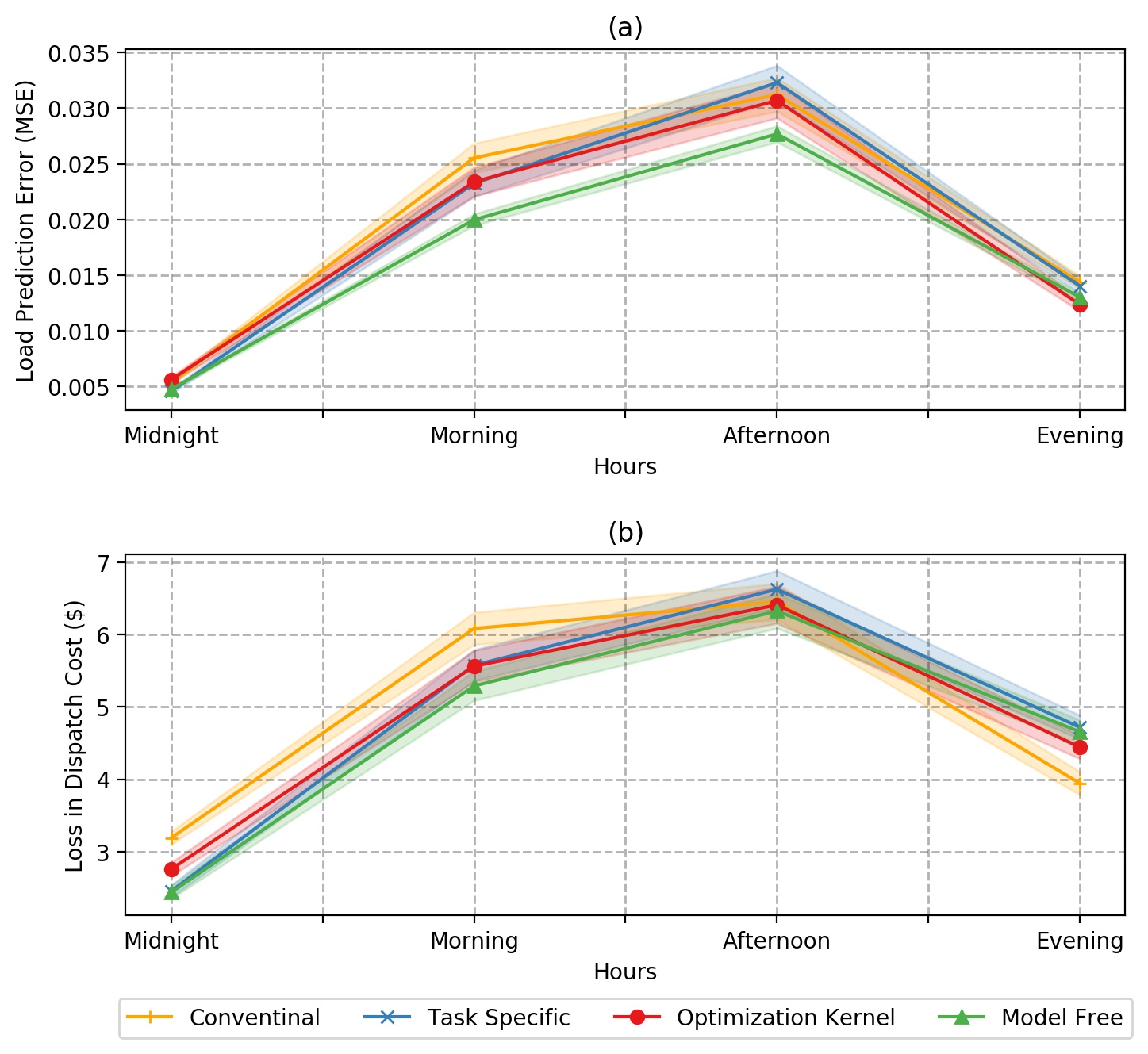}}
\caption{Effectiveness Comparison (the shadowed areas illustrate the $\pm0.1\sigma$ zone for the four frameworks).}
\label{difm}
\end{figure}

During the comparison, we divide each day into four periods: midnight (from 0:00 am to 6:00 am), morning (from 6:00 am to 12:00 pm), afternoon (from 12:00 pm to 6:00 pm), and evening (from 6:00 pm to 0:00 am), and conduct the comparison for these four periods, respectively. Figure \ref{difm} plots the comparison with respect to the two evaluation metrics.

It is interesting to observe that the conventional approach doesn't produce the best prediction in all periods. Yet it does incur relatively high loss in the dispatch cost (i.e., additional dispatch cost). As suggested by Fig.\ref{difm}(b), all the other three frameworks improve the effectiveness of the conventional approach (benchmark). More precisely, task specific criteria achieves $1.53 \%$ less loss than the benchmark; optimization kernel framework further improves the effectiveness by $2.13\%$; and the model-free end-to-end learning framework achieves $1.27\%$ additional improvement by relaxing the assumption of the load prediction distributions. This highlights the remarkable performance of our proposed frameworks.




\subsection{Performance Evaluation: Robustness}
\begin{figure}[t]
\centerline{\includegraphics[width=8.5cm]{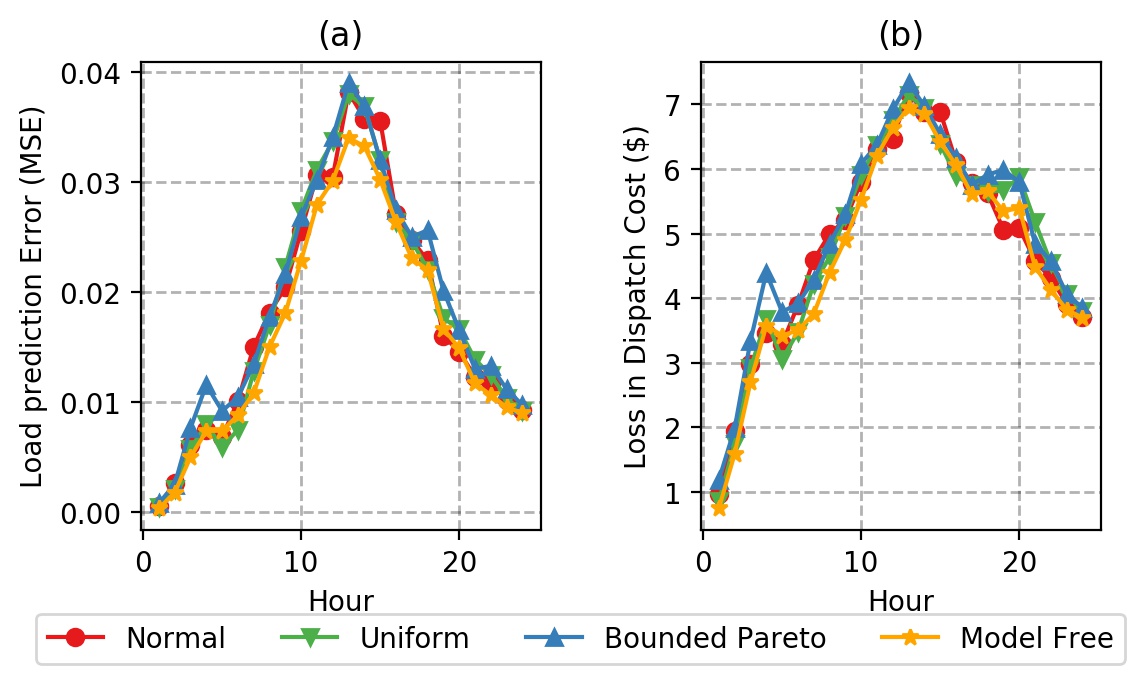}}
\caption{Robustness Evaluation for End-to-End Framework with Optimization Kernel.}
\label{robust}
\end{figure}

The task specific framework and the optimization kernel framework both requires the knowledge on the type of prediction distribution for effective learning. We target to examine their robustness through synthetic data generated from Normal distribution, Uniform distribution, and Bounded Pareto distribution. We use the model-free framework as the benchmark.
We take the optimization kernel framework as an example. Fig. \ref{robust} visualizes its robustness evaluation. It is suggested that the proposed framework is rather robust to most light tail distributions (e.g., Normal distribution and Uniform distribution) as the loss/error in the two evaluation metrics is bounded by $4 \%$ compared with the benchmark. The robustness is weakened facing heavy tail distributions (e.g., bounded Pareto distribution) yet the loss/error is still bounded by $8 \%$.



\subsection{Performance Evaluation: Efficiency}
\begin{figure}[t]
\centerline{\includegraphics[width=5cm]{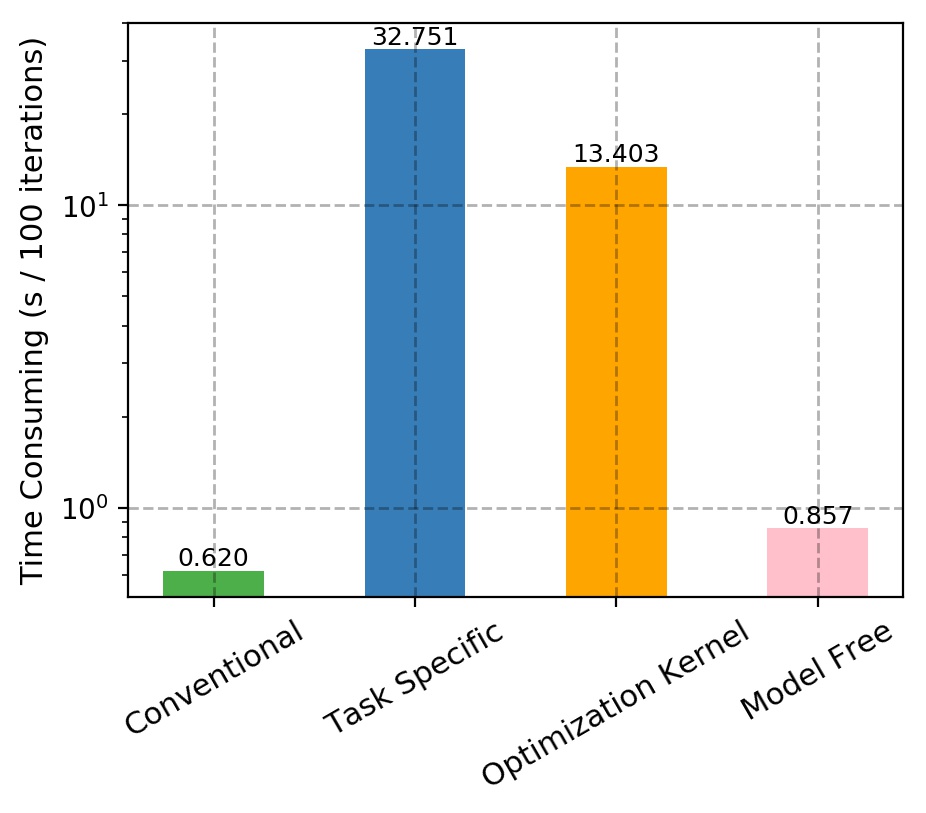}}
\caption{Efficiency Comparison.}
\label{time}
\end{figure}

While the three proposed frameworks all outperform the conventional approach, the extra performance does not come free: the learning process is more complex and hence more time consuming than that of the conventional approach. Figure \ref{time} illustrates the runing time for the four frameworks on a log-scale. Compared with the task specific framework, the optimization kernel speeds up the learning process by $182\%$, which verifies the effectiveness of the kernel. The model-free framework further speeds up the process to the comparable level of the conventional approach!

\subsection{Performance Evaluation for IEEE 39-bus System}
We further verify the performance of the model-free framework on a larger system: the IEEE 39-bus system. We follow the network parameters in the system as suggested in \cite{39bus}. We set the marginal cost of the generators from bus 30 to bus 39 in an increasing order: from $\$30$/MWh up to $\$48$/MWh at a step size of $\$2$/MWh. As such, we set unit shortage penalty $\gamma_1$ to be $\$50$/MWh and the unit excess penalty to be $\$2$/MWh.

\begin{figure}[t]
\centerline{\includegraphics[width=8.5cm]{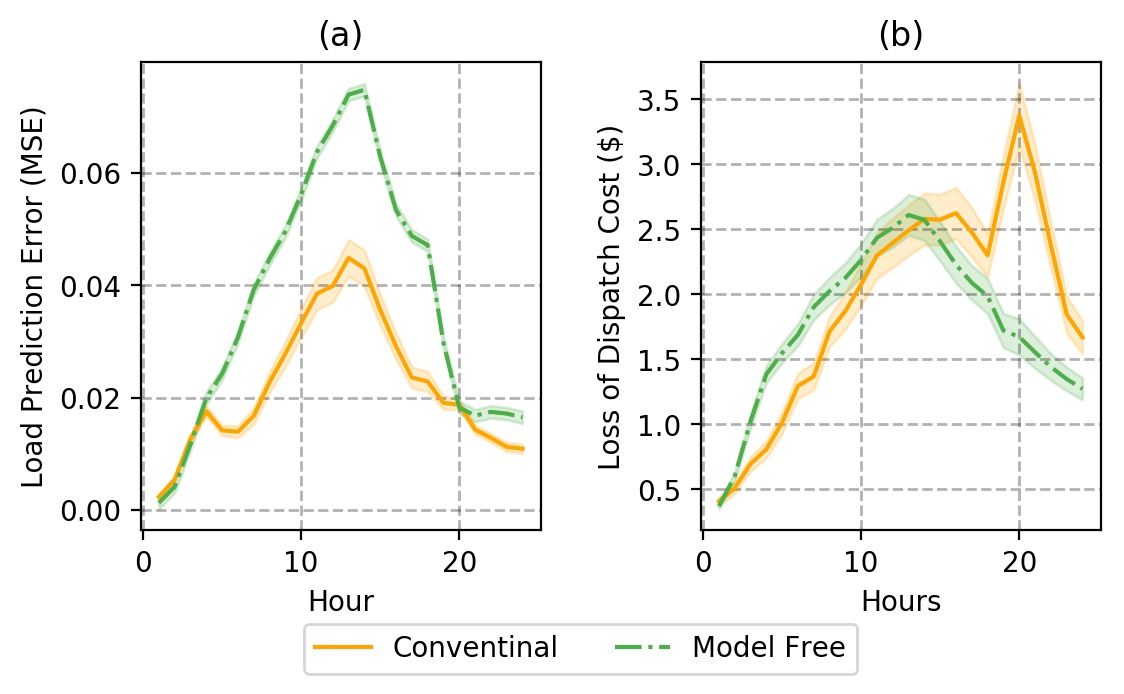}}
\caption{Performance Evaluation for IEEE 39-Bus System  (the shadowed areas illustrate the $\pm0.1\sigma$ zone).}
\label{IEEE39}
\end{figure}

Figure \ref{IEEE39} illustrates the performance improvement of the model-free framework, compared with the conventional approach as a benchmark. While our framework generates $72.85\%$ more load prediction error, it does reduce the loss in dispatch cost by $8.24\%$.




\section{Conclusion}
\label{conclude}
In this paper, we seek to design an effective end-to-end machine learning framework for ED. Specifically, motivated by the task specific criteria, we design the optimization kernel to speed up the learning process, which ultimately leads to our model-free framework. Numerical studies verify the efficiency and effectiveness of our proposed schemes.

Our work can be extended in various interesting ways. For example, we consider the network constraints satisfies on expectation. That is, we implicitly assume the prediction error can be bounded within a relatively small range, which won't affect the network constraints too much. It will be interesting to propose a more adaptive way to ensure the network constraints in the end-to-end learning framework.



\begin{appendices} 
\section{} 
\subsection{Proof for Lemma 2}
First, we prove the funtion domain of $g_t$ is $[g_t^{\max}$, $g_t^{\min}]$.

Define $$\boldsymbol{g}_t^{\min} = (g_{1t}^{\min}, g_{2t}^{\min}, ... , g_{nt}^{\min}),$$ $$\boldsymbol{g}_t^{\max} = (g_{1t}^{\max}, g_{2t}^{\max}, ... , g_{nt}^{\max}).$$ These two dispatch policies correspond to $g_t^{\min}$ and $g_t^{\max}$, respectively.

Note that, any $g_t^\alpha \in (g_t^{\min}, g_t^{\max})$ can be expressed as follows:
\begin{align}
    g_t^\alpha = (1-\alpha)g_t^{\min} + \alpha g_t^{\max}\quad(0<\alpha<1)
\end{align}

We can construct the following generation policy $g_t^\alpha$:
\begin{align}
    \boldsymbol{g}_t^{\alpha} = (1-\alpha)\boldsymbol{g}_t^{\min} + \alpha \boldsymbol{g}_t^{\max}.
\end{align}

Since constraints (\ref{eq2}) and (\ref{eq3}) are both linear, it is straightforward to verify that $\boldsymbol{g}_t^\alpha$ is a feasible solution to $C_t(g_t)$. Moreover, the arbitrariness of $\alpha$ implies that the range of $g_t$ is $(g_t^{\min},g_t^{\max})$.



The continuity and piece-wise linearity of $C_t(g_t)$ are immediate results of \cite{holder}.
Then we prove the convexity of $C_t(g_t)$:

Suppose $\boldsymbol{g}_t^a$ and $\boldsymbol{g}_t^b$ are optimal dispatch policies of ${g}_t^a$ and ${g}_t^b$. Without loss of generality, assume:
\begin{align}
    {g}_t^a < {g}_t^b.
\end{align}
We can construct a policy $\boldsymbol{g}_t^\beta$ from $\boldsymbol{g}_t^a$ and $\boldsymbol{g}_t^b$ as follows:
\begin{align}
    \boldsymbol{g}_t^\beta = (1-\beta)&\boldsymbol{g}_t^a + \beta\boldsymbol{g}_t^b\quad(0 < \beta < 1)\\
    &g_t^a < g_t^\beta < g_t^b
\end{align}
The linearity in $C_t(\cdot)$ further implies that:
\begin{align}
    C(\boldsymbol{g}_t^\beta) = (1-\beta)C(&\boldsymbol{g}_t^a) + \beta{C}(\boldsymbol{g}_t^b)
\end{align}
Hence,
\begin{equation}
    \begin{split}
    C_t(g_t^\beta) &\leq C(\boldsymbol{g}_t^\beta)\\
&= (1-\beta)C(\boldsymbol{g}_t^a) + \beta{C}(\boldsymbol{g}_t^b)\\
&= (1-\beta)C_t({g}_t^a) + \beta{C_t}({g}_t^b).   
    \end{split}
\end{equation}

The first inequality is due to the definition of $C_t(g_t)$. This concludes our proof.
\end{appendices} 
$\hfill\blacksquare$ 
\subsection{Proof for Lemma 3}
As $g_t^{k}$ and $g_t^{k+1}$ are two adjacent breaking points, we know:
\begin{align}
    C_t({g}_t^\gamma) = (1-\gamma)C_t({g}_t^{k}) + \gamma{C_t({g}_t^{k})}.
\end{align}
On the other hand, the linearity in the cost function $C(\cdot)$ implies:
\begin{equation}
    \begin{split}
        C(\boldsymbol{g}_t^\gamma) &= (1-\gamma)C(\boldsymbol{g}_t^{k}) + \gamma{C(\boldsymbol{g}_t^{k})}\\
    &= (1-\gamma)C_t({g}_t^{k}) + \gamma C_t({g}_t^{k})\\
    &= C_t(g_t^\gamma).
    \end{split}
\end{equation}
$\hfill\blacksquare$ 
\subsection{Proof for Lemma 4}
In Lemma 2, we prove that $C_t(g_t)$ is convex, which means:
\begin{align}
    C_t''(x) \geq 0.
\end{align}
Also, $H(x)$ is a \emph{cdf}, which implies its derivative is non negative. Hence,
\begin{align}
    K'(x) = (\gamma_1 + \gamma_2)H'(x) + C_t''(x)\geq 0.
\end{align}
When $g_t^{\min} < g^* < g_t^{\max}$, we know:
\begin{align}
    K(g^*) = 0.
\end{align}
Hence, $K(x)$ is either always positive over $[g^* , g(\hat{\boldsymbol{\theta}}_t)]$ or always negative over $[g(\hat{\boldsymbol{\theta}}_t), g^*]$.

Following the same routine, we can verify this argument when $g^* = g_t^{\min}$ and $g^* = g_t^{\max}$.$\hfill\blacksquare$ 
\end{document}